\documentclass{article}
\usepackage[utf8]{inputenc}
\usepackage{mathtools}
\usepackage{amsmath,amssymb,amsbsy}

\newcommand{\zet}{\mathbb{Z}}

\newcommand{\qed}{\hfill \rule{.1in}{.1in}}

\newtheorem{thm}{Theorem}[section]
\newtheorem{lem}[thm]{Lemma}

\newtheorem{cor}[thm]{Corollary}

\newtheorem{dfn}[thm]{Definition}

\newtheorem{conj}[thm]{Conjecture}

\begin{document}
\title{Irregular labeling on Abelian groups of digraphs}

\author{Sylwia Cichacz\footnote{This work was partially supported by the Faculty of Applied Mathematics AGH UST statutory tasks within subsidy of Ministry of Science and Higher Education.}
\normalsize \\AGH University of Science and Technology, \vspace{2mm} Poland\\
}

\maketitle
\begin{abstract}
Let  $\overrightarrow{G}$  be  a  directed  graph  of order $n$ with no component of order
 less than $4$, and let $\Gamma$ be a finite Abelian group such that
  $|\Gamma|\geq n+6$.
We  show that  there exists a
 mapping $\psi$ from the arc set $E(\overrightarrow{G})$ of
$\overrightarrow{G}$ to an Abelian group $\Gamma$
such
that if we define a mapping $\varphi_{\psi}$ from the vertex set $V(\overrightarrow{G})$ of
$\overrightarrow{G}$
to $\Gamma$
by
$$\varphi_{\psi}(x)=\sum_{y\in N^+(x)}\psi(xy)-\sum_{y\in N^-(x)}\psi(yx),\;\;\;(x\in V(\overrightarrow{G})),$$
then $\varphi_{\psi}$
is injective. Such a labeling $\psi$ is called \textit{irregular}.
\end{abstract}

 \section{Introduction}

Let  $\overrightarrow{G}=(V,E)$ be a directed graph without isolated vertices.
For arcs we omit the arrow on the top from the notation. An arc $xy$
 is considered to be directed from $x$ to $y$, moreover $y$ is called the \emph{head} and $x$ is called the \emph{tail} of the arc. For a vertex $x$, the set of head endpoints adjacent to $x$ is denoted by $N^+(x)$, and the set of tail endpoints adjacent to $x$ is  denoted by $N^-(x)$.

Assume $\Gamma$ is a finite Abelian group with the operation denoted by~$+$.
  For convenience
we will write $ka$ to denote $a + a + \ldots + a$ where the element $a$ appears $k$ times, $-a$ to denote the inverse of $a$, and
we will use $a - b$ instead of $a+(-b)$.  Moreover, the notation $\sum_{a\in S}{a}$ will be used as a short form for $a_1+a_2+a_3+\dots$, where $a_1, a_2, a_3, \dots$ are all elements of the set $S$. The identity element of $\Gamma$ will be denoted by $0$. 

Suppose that there exists a
 mapping $\psi$ from the arc set $E(\overrightarrow{G})$ of
$\overrightarrow{G}$ to an Abelian group $\Gamma$
such
that if we define a mapping $\varphi_{\psi}$ from the vertex set $V(\overrightarrow{G})$ of
$\overrightarrow{G}$
to $\Gamma$
by
$$\varphi_{\psi}(x)=\sum_{y\in N^+(x)}\psi(xy)-\sum_{y\in N^-(x)}\psi(yx),\;\;\;(x\in V(\overrightarrow{G})),$$
then $\varphi_{\psi}$
is injective. In this situation, we say that $\overrightarrow{G}$ is \textit{realizable} in
$\Gamma$, and that the mapping $\psi$ is $\Gamma$-\textit{irregular}.

The corresponding problem in the case of simple graphs was considered
 in \cite{ref_AnhCic,ref_AnhCicPrzy}.
  For  $\Gamma=(\zet_2)^m$ the problem was raised in \cite{Tuza}.
We easily see that if $\overrightarrow{G}$ is realizable in $(\zet_2)^m$,
 then every component of $\overrightarrow{G}$ has order at least 3
 (recall that we are assuming $\overrightarrow{G}$ has no isolated vertex).
The following results have been shown:

\begin{thm}[\cite{ref_CaccJia,Egawa}] \label{Egawa} Let\/ $\overrightarrow{G}$ be a directed
 graph of order\/ $n$, with no weakly component of order less than\/ $3$.
Then\/ $\overrightarrow{G}$ is realizable in\/ $(\mathbb{Z}_2)^m$ if and only if\/
  $n\leq 2^m$ and\/ $n\neq 2^m-2$. 
\end{thm}

We note that the case of $n=2^m-1$ means realization in
 $(\mathbb{Z}_2)^m\setminus \{0\}$ because the sum of all elements
 in $(\mathbb{Z}_2)^m$ is zero for $m\geq 2$, and consequently the single
  element missing from the domain of $\varphi_{\psi}$ must be 0.

\begin{thm}[\cite{Fukuchi}]  Let\/ $p$ be an odd prime and let\/
 $m\geq 1$ be an integer.
If\/~$\overrightarrow{G}$ is a directed graph of order\/ $n$ without isolated vertices such that\/
 $n \leq p^m$, then\/  $\overrightarrow{G}$ is realizable
  in\/ $(\zet_p)^m$.
\end{thm}

Cichacz and Tuza showed that if  $n$ is large enough
with respect to an arbitrarily fixed $\varepsilon > 0$ then $\overrightarrow{G}$
has a $\Gamma$-irregular labeling for any $\Gamma$ such that $|\Gamma|>(1+\varepsilon)n$ \cite{CicTuz}.  In the case of any graph they obtain the following:
\begin{thm}[\cite{CicTuz}]\label{glowne}Any directed graph\/ $\overrightarrow{G}$ of order\/ $n$ with no weakly connected components
 of cardinality less than\/ $3$ has a $\Gamma$-irregular labeling for every\/ $\Gamma$ such that\/
  $|\Gamma|\geq 2n+2\sqrt{n-1/2}-1$.
\end{thm}

Recently Cichacz and Suchan showed that if all weakly components are bigger that $3$, then the coefficient $2\sqrt{n-1/2}-1$ can be improved:
\begin{thm}[\cite{CicrSuch2}]Any directed  graph\/ $\overrightarrow{G}$ of order\/ $n$ with no weakly connected components
 of order less than\/ $4$ has a $\Gamma$-irregular labeling for every\/ $\Gamma$ such that\/
  $|\Gamma|\geq 2n+1$.
\end{thm}
Moreover for groups with $|I(\Gamma)|\neq 1$ the authors showed better bound. 
\begin{cor}[\cite{CicrSuch2}]\label{bez} Any digraph\/ $\overrightarrow{G}$ of order\/ $n$ with no weakly connected components
 of cardinality less than\/ $4$ has a $\Gamma$-irregular labeling for every\/ $\Gamma$ such that\/
  $|I(\Gamma)|\neq 1$ and $|\Gamma|\geq n+5$.\label{liniowe}\end{cor}

They also stated the following conjecture:

\begin{conj}[\cite{CicrSuch2}]There exists a constant $K$ such that any digraph\/ $\overrightarrow{G}$ of order\/ $n$ with no weakly connected components
 of order less than\/ $3$ has a $\Gamma$-irregular labeling for every\/ $\Gamma$ such that $|\Gamma|\geq n+K$.
\end{conj}

In this paper we show that the conjecture is true for any digraph with no weakly connected components
 of order less than\/ $4$, namely we prove  any directed  graph\/ $\overrightarrow{G}$ of order\/ $n$ with no weakly connected components
 of order less than\/ $4$ has a $\Gamma$-irregular labeling for every\/ $\Gamma$ such that\/
  $|\Gamma|\geq n+6$.
  
  \section{Preliminaries}

 A non-trivial
finite group has elements of order $2$ (involutions) if and only if the order of the group is even. The fundamental theorem of finite Abelian groups states that a finite Abelian
group $\Gamma$ of order $n$ can be expressed as the direct product of cyclic subgroups of prime-power order. This implies that
$$\Gamma\cong\zet_{p_1^{\alpha_1}}\times\zet_{p_2^{\alpha_2}}\times\ldots\times\zet_{p_k^{\alpha_k}}\;\;\; \mathrm{where}\;\;\; n = p_1^{\alpha_1}\cdot p_2^{\alpha_2}\cdot\ldots\cdot p_k^{\alpha_k}$$
and $p_i$ for $i \in \{1, 2,\ldots,k\}$ are not necessarily distinct primes. This product is unique up to the order of the direct product.  Since the properties and results in this paper are invariant under the isomorphism ($\cong$) between groups, we only need to consider one group in each isomorphism class.

Recall that any group element $\iota\in\Gamma$ of order 2 (i.e., $\iota\neq 0$ and $2\iota=0$) is called an \emph{involution}. Let us denote the number of involutions in $\Gamma$ by $|I(\Gamma)|$. By the fundamental theorem of Abelian groups, it is easy to see that any Abelian group $\Gamma$ can be factorized as $\Gamma\cong L\times H$,  with $|L|=2^\eta$ for a nonnegative integer $\eta$ and $|H|=\rho$ for an odd positive integer $\rho$. Note that in this case $|I(\Gamma)|=|I(L)|$. Thus if $|I(\Gamma)|=1$ then $\Gamma\cong \zet_{2m}\times H$ for $m\geq 1$ and $|H|$  odd.

Recall that the sum of all elements of a group $\Gamma$ is equal to the sum of its involutions and the identity element. The following lemma was proved in~\cite{ref_ComNelPal} (see \cite{ref_ComNelPal}, Lemma 8).

\begin{lem}[\cite{ref_ComNelPal}]\label{involutions} Let  Let $\Gamma$  be an Abelian group.
\begin{itemize}
 \item[-] If $\Gamma$ has exactly one involution $\iota$, then $\sum_{g\in \Gamma}g= \iota$.
\item[-] If $\Gamma$ has no involutions, or more than one involution, then $\sum_{g\in \Gamma}g=0$.
\end{itemize}
\end{lem}

\section{Zero-sum partition}
A subset $S$ of $\Gamma$ is called a zero-sum subset if $\sum_{a\in S} a=0$.
It turns out that a realization of $\overrightarrow{G}$ in an Abelian group $\Gamma$
 is strongly connected with a zero-sum partition of $\Gamma$ \cite{ref_AigTri,CicTuz,Fukuchi}. It was shown the following.
\begin{lem}[\cite{CicTuz}]\label{zerosum} A directed graph\/ $\overrightarrow{G}$ with no isolated
 vertices has a $\Gamma$-irregular labeling if and only if there exists an injective
mapping\/ $\varphi$ from $V(\overrightarrow{G})$ to $\Gamma$ such that\/
 $\sum_{x\in C}\varphi(x)=0$ for every  weakly connected component $C$ of\/ $\overrightarrow{G}$.
\end{lem}
Since we will often consider subsets of fixed cardinalities throughout the paper, let us present an abbreviated notation. Given a set, any of its subsets of cardinality $k$ is called a \textit{$k$-subset}.
\begin{cor}[\cite{Tannenbaum1,Zeng}]\label{Zeng}
{Let} $n$ be even positive natural number and  $m,l$ be natural numbers such that $3m+2l = n-2$. Then the set $S=\zet_n\setminus\{0,\frac{n}{2}\}$ can be partitioned into $m$ zero-sum $3$-sets $A_1,A_2,\ldots,A_{m}$ and $l$ zero-sum $2$-sets $B_1,B_2,\ldots,B_l$.  
\end{cor}

We call a $6$-subset $C$ of an Abelian group $\Gamma$ \textit{good} if $C = \{c, d,-c -d,-c,-d, c + d\}$ for some $c$ and $d$ in $\Gamma$. Notice that the sum of elements of a good $6$-subset is $0$. Moreover, it can be partitioned into three zero-sum $2$-subsets or two zero-sum $3$-subsets.

 The following definition was given by Tannenbaum \cite{Tannenbaum1}.

\begin{dfn}\label{dfn:skolem}Let $\Gamma$ be a finite Abelian group of order $m=6k+s$ for a non-negative integer $k$ and $s\in\{1,3,5\}$. A partition of $\Gamma\setminus\{0\}$ into $k$ good $6$-subsets and $(s-1)/2$ zero-sum $2$-subsets is called a {\em Skolem partition of $\Gamma$}.
\end{dfn}

Tannenbaum showed the following:

\begin{thm}[\cite{Tannenbaum1}]\label{Tannenbaum1}Let $\Gamma$ be a finite Abelian group such that $|I(\Gamma)|=0$, then $\Gamma$ has a Skolem partition. 
\end{thm}

Let us introduce some notation.  Let $L$ and $H$ be Abelian groups such that $L$ is cyclic group  with exactly one involution $\iota$ and $|H|$ is odd.
 Observe that $|L|=2l+2+ 6m$ and $|H|=2p+1+ 6k$ for some integers $m,k$ and $l,p\in\{0,1,2\}$. By Corollary~\ref{Zeng} the set $S=L\setminus\{0,\iota\}$ can be partitioned into $m$ zero-sum $3$-subsets $A_1,A_2,\ldots,A_{m}$ and $l$ zero-sum $2$-subsets $B_1,B_2,\ldots,B_l$.  Denote $A_i=\{a_{i,0}, a_{i,1},a_{i,2}\}$ for $i=1,2,\ldots,m$ and $B_i=\{b_i,-b_i\}$ for $i=1,2,\ldots,l$. By Theorem~\ref{Tannenbaum1}, there exists a Skolem partition of $H$, hence $H\setminus\{0\}$ can be partitioned into $k$ good $6$-subsets $C_1,C_2,\ldots,C_{k}$ and $p$ zero-sum $2$-subsets $D_1,D_2,\ldots,D_p$.   Denote $$C_i=\{c_{j,0}, c_{j,1},-c_{j,0} - c_{j,1},-c_{j,0},-c_{j,1},c_{j,0} + c_{j,1}\}$$ for $i=1,2,\ldots,k$ and $D_i=\{d_i,-d_i\}$ for $i=1,\ldots,p$.
 Let $c_{i,2}=-(c_{i,0}+c_{i,1})$ for $i=1,2,\ldots,k$.
Let now
$$W_0=\bigcup\limits_{i=1}^{m}\bigcup\limits_{h=0}^{2}\bigcup\limits_{j=1}^{k}\{(a_{i,1},c_{j,h}),(a_{i,2}, c_{j,1+h}),(a_{i,3}, c_{j,2+h}),
(-a_{i,1},-c_{j,h}),(-a_{i,2}, -c_{j,1+h}),(-a_{i,3}, -c_{j,2+h})\}$$ 
$$W_1=\bigcup\limits_{j=1}^{k}\{(0,c_{j,1}),(0, c_{j,2}),(0,c_{j,3}),(0,-c_{j,1}),(0, -c_{j,2}),(0,-c_{j,3})\},$$
$$W_2=\bigcup\limits_{h=0}^{2}\bigcup\limits_{j=1}^{k}\{(b_1,c_{j,h}),(-b_1, c_{j,1+h}),(0,c_{j,2+h}),(-b_1,-c_{j,h}),(b_1, -c_{j,1+h}),(0,-c_{j,2+h})\},$$

$$T=\bigcup\limits_{h=0}^{2}\bigcup\limits_{i=1}^{m}\{(a_{i,h},d_1),(a_{i,1+h},-d_1),(a_{i,2+h},0)\},$$
where the second subscripts are taken modulo 3.

Note that $|W_0|=(|L|-2l-2)(|H|-2p-1)$, $|W_1|=|H|-2p-1$, $|W_2|=3(|H|-2p-1)$ for $l\neq 0$ and $|W_2|=0$ otherwise, $|T|=3(|L|-2l-2)$ for $p\neq 0$ and $|T|=0$ otherwise.

The main result of this section is the following.

\begin{thm}   \label{zero-sum}

Let $\Gamma$ be of order\/ $n$, $|I(\Gamma)|=1$, $\iota_{\Gamma}$ be the involution in $\Gamma$ and consider any integers\/
 $r_1,r_2,\dots,r_t$ with\/ $n -2=r_1 + r_2 + \ldots + r_t$
  and {$r_i \geq 4$} for all\/ $1 \leq i \leq t$.
Then
 there exist pairwise disjoint zero-sum subsets\/ $A_1, A_2,\ldots , A_t$ in\/ $\Gamma\setminus\{0,\iota_{\Gamma}\}$
 such that\/ $|A_i| = r_i$ for all\/ $1 \leq i \leq t$.
\end{thm}
\textit{Proof. } Note that $\Gamma\cong L\times H$, where $L\cong \zet_{2^{\eta}}$ for $\eta\geq 1$ and $|H|$ is odd. By Corollary~\ref{Zeng} we may assume that $|H|= 9$ or $|H|\geq 25$ (because otherwise $\Gamma$ is cyclic). Let $\iota$ be the involution in $L$ (then $(\iota,0)=\iota_{\Gamma}$).

 Assume that $r_1,\dots,r_s$ are all odd, and $r_{s+1},\dots,r_t$ are all even. Then $s$ is even because $n$ is even and $s\leq \left\lfloor\frac{|\Gamma|-2}{5}\right\rfloor$.

\textit{Case 1.} $|L|>2$

For $|L|=4$ (i.e. $l=1$) there is
 $\frac{|W_2|}{3}=|H|-2p-1\geq \left\lfloor\frac{4|H|-2}{5}\right\rfloor=\left\lfloor\frac{|\Gamma|-2}{5}\right\rfloor\geq s$ for $|H|\geq 9$. Thus we can find $s$ disjoint zero-sum $3$-sets $T_1,\ldots,T_s$ in $W_2$, the set $\Gamma\setminus\left(\{(0,0),(\iota,0)\}\cup \bigcup\limits_{i=1}^{s}T_i\right)$ contains zero-sum $2$-sets, thus there exist pairwise disjoint zero-sum subsets\/ $A_1, A_2,\ldots , A_t$ in\/ $\Gamma\setminus\{(0,0),(\iota,0)\}$
 such that\/ $|A_i| = r_i$ for all\/ $1 \leq i \leq t$.
 
 From now on we assume that $|L|>4$, what implies that $W_0\neq \emptyset$. Let $W=W_0\cup W_1$ for  $|L|\equiv 2\pmod6$ (i.e. $l=0$) and  $W=W_0\cup W_2$ for  $|L|\equiv 4\pmod6$ (i.e $l=1$). Assume first that $|H|\equiv 1 \pmod 6$ (i.e. $p=0$), then  one can easily check  that $\frac{(|H|-1)(|L|-1)}3=\frac{|W|}3\geq    \frac{|H||L|-2}{5}=\frac{|\Gamma|-2}{5}\geq s$ since $|L|\geq 8$ and $|H|\geq 3$. Therefore can find $s$ disjoint zero-sum $3$-sets $T_1,\ldots,T_s$ in $W$, the set $\Gamma\setminus\left(\{(0,0),(\iota,0)\}\cup \bigcup\limits_{i=1}^{s}T_i\right)$ contains zero-sum $2$-sets, thus there exist pairwise disjoint zero-sum subsets\/ $A_1, A_2,\ldots , A_t$ in\/ $\Gamma\setminus\{(0,0),(\iota,0)\}$
 such that\/ $|A_i| = r_i$ for all\/ $1 \leq i \leq t$.

 Assume now that $|H|\not\equiv 1\pmod 6$ (i.e $p\in\{1,2\}$).  If now  $s\leq |W|/3$, then we can find $s$ disjoint zero-sum $3$-sets $T_1,\ldots,T_s$ in $W$, the set $\Gamma\setminus\left(\{(0,0),(\iota,0)\}\cup \bigcup\limits_{i=1}^{s}T_i\right)$ contains zero-sum $2$-sets, thus there exist pairwise disjoint zero-sum subsets\/ $A_1, A_2,\ldots , A_t$ in\/ $\Gamma\setminus\{(0,0),(\iota,0)\}$
 such that\/ $|A_i| = r_i$ for all\/ $1 \leq i \leq t$.
 
 From now $s>|W|/3$. We find $(|L|-2l-2)=\frac{|T|}{3}$ zero-sum $3$-subsets in $T$ because $|T|\leq |W|$ for $|H|\geq 9$. Let $s'=s-m$.
 Note that $\frac{|W|+|T|}{3}=\frac{(|H|-2p+2)(|L|-1)-3-6l}{3}\geq \frac{|L||H|-2}{5}\geq s$ for $|H|>5$ and $|L|>4$. The remaining $s-(|L|-2l-2)$ zero-sum $3$-sets we find  in $W$. Hence all the other elements in $\Gamma\setminus\{(0,0),(\iota,0)\}$ form zero-sum $2$-subsets,  we are done. 
 
  \textit{Case 2.} $|L|=2$.
 
  \textit{Case 2.1} $|H|\neq3^{\beta}$ for any positive integer $\beta\geq 2$. 
  
  Suppose first that $\Gamma \cong\zet_2\times \zet_5\times \zet_5$, observe that in that case $s\leq 8$. The set $H\setminus\{0\}=\zet_5\times \zet_5$  is a union of 4 zero-sum good $6$-subsets by Theorem~\ref{Tannenbaum1}. Thus for $L=\zet_2$ the set $W_1$ has also  4 zero-sum good $6$-subsets and we are done analogously as in Case 1.
  
  Therefore there exists $n_1\not\equiv3\pmod 6$, $n_1\geq5$ such that $H=\zet_{n_1}\times H'$ for $|H'|\geq 7$ (otherwise the group $\Gamma$ is cyclic and we are done by Lemma~\ref{Zeng}). 
 Let $L'=L\times \zet_{n_1}$, then $|L'|\geq 10$ and  $|L'|=2l'+2+ 6m'$ and $|H'|=2p'+1+ 6k'$ for some integers $m',k'$ and $l'\in\{1,2\}$ and $p'\in\{0,1,2\}$. We are done as in Case 1.

\textit{Case 2.2} $|H|=3^{\beta}$ for some positive integer $\beta\geq 2$. 
If  $\Gamma\cong \zet_2\times H$ for $|H|=9$, then $s\in\{0,2\}$. Since $H=\{c_1,c_2,-c_1-c_2,-c_1,-c_2,c_1+c_2\}\cup\{0,d,-d\}$ by Theorem~\ref{Tannenbaum1} the set $S=\{(0,c_1),(0,c_2),(0,-c_1-c_2),(0,-c_1),(0,-c_2),(0,c_1+c_2)\}\subset\Gamma\setminus(\{(0,0),(\iota,0)\}$ is a good 6-subset  and $\Gamma\setminus(\{(0,0),(\iota,0)\}\cup S)$ form zero-sum $2$-subsets,  we are done. Hence we assume that $|H|\geq 27$.
 
 Since  $H$ is not cyclic we have   $H\cong \zet_{n_1}\times H'$ for $n_1=3^{\eta}$ for some $0<\eta<\beta$ and  $|H'|=3^{\beta-\eta}$. 
By Theorem~\ref{Tannenbaum1}, there exists a Skolem partition of $H'$, hence $H'\setminus\{0\}$ can be partitioned into $k'=(|H'|-3)/6$ good $6$-subsets $C_1',C_2',\ldots,C_{k'}'$ and one zero-sum $2$-subset $D'_1=\{-d'_1,d'_1\}$.   Denote $C_i'=\{c_{j,0}', c_{j,1}',c_{j,2}',-c_{j,0}
,-c_{j,1}',-c_{j,2}'\}$ for $i=1,2,\ldots,k'$, where $c_{j,2}'=-c_{j,0}'-c_{j,1}'$. Observe that  for $L'\cong \zet_2\times\zet_{n_1}\cong\zet_{2n_1}$ there is $|L'|=6+ 6m'$ for $m'=(n_1-3)/3$. By Corollary~\ref{Zeng} the set $L'\setminus\{0,\iota'\}$ (where $\iota'$ is the involution in $L'$) can be partitioned into $m'$ zero-sum $3$-subsets $A_1',A_2',\ldots,A_{m'}'$ and $2$ zero-sum $2$-subsets $B_1,B_2$.  Denote $A_i'=\{a_{i,0}', a_{i,1}',a_{i,2}'\}$ for $i=1,2,\ldots,m'$ and $B_i'=\{b_i',-b_i'\}$ for $i=1,2$.

Suppose first  that $n_1=3$.  Note that $\Gamma\cong L'\times H'$ for  $L'\cong \zet_6=\{0,1,2,3,4,5\}$ and $|H'|=|H|/3\geq 9$. 

Let 
$$W_1'=\bigcup\limits_{h=0}^{2}\bigcup\limits_{j=1}^{k'}\{(1,c'_{j,h}),(5, c'_{j,1+h}),(0,c'_{j,2+h}),(5,-c'_{j,h}),(1, -c'_{j,1+h}),(0,-c'_{j,2+h})\},$$

$$W_3'=\bigcup\limits_{j=1}^{k'}\{(2,c'_{j,1}),(2, c'_{j,2}),(2,c'_{j,3}),(4,-c'_{j,1}),(4, -c'_{j,2}),(4,-c'_{j,3})\},$$
 $$S'=\{(4,d),(4,0),(4,-d),(2,-d),(2,0),(2,d)\}.$$
Note that $W'=W_1'\cup W_3'\cup S'\subset\Gamma \setminus(\{(0,0),(\iota',0)\}$ are good $6$-subsets and $\frac{|W'|}{3}=\frac{5(|H'|-3)+6}{3}\geq \frac{3|H'|-2}{5}\geq s$ for $|H'|\geq 9$. Therefore can find $s$ disjoint zero-sum $3$-sets $T_1,\ldots,T_s$ in $W'$, the set $\Gamma\setminus\left(\{(0,0),(\iota,0)\}\cup \bigcup\limits_{i=1}^{s}T_i\right)$ contains zero-sum $2$-sets, thus there exist pairwise disjoint zero-sum subsets\/ $A_1, A_2,\ldots , A_t$ in\/ $\Gamma\setminus\{(0,0),(\iota',0)\}$
 such that\/ $|A_i| = r_i$ for all\/ $1 \leq i \leq t$.

 From now on $n_1\geq 9$. Define sets $W_0$, $W_2$ and $T$ as in Case 1 (for $L'$ and $H'$ instead of $L$ and $H$). Let $W=W_0\cup W_2$, then $|W|=(|L'|-3)(|H'|-3)$. If now  $s\leq |W|/3$, then we can find $s$ disjoint zero-sum $3$-sets $T_1,\ldots,T_s$ in $W$, the set $\Gamma\setminus\left(\{0,\iota_{\Gamma}\}\cup \bigcup\limits_{i=1}^{s}T_i\right)$ contains zero-sum $2$-sets, thus there exist pairwise disjoint zero-sum subsets\/ $A_1, A_2,\ldots , A_t$ in\/ $\Gamma\setminus\{(0,0),(\iota,0)\}$
 such that\/ $|A_i| = r_i$ for all\/ $1 \leq i \leq t$.
 
 From now $s>|W|/3$. We find $(|L'|-6)=\frac{|T|}{3}$ zero-sum $3$-subsets in $T$ because $|T|\leq |W|$ for $|H'|\geq 9$. Let $s'=s-m$.
 Note that $\frac{|W|+|T|}{3}=\frac{|H'|(|L'|-3)-9}{3}\geq \frac{|L'||H'|-2}{5}\geq s$ for $|H'|\geq 9$ and $|L'|\geq 18$. The remaining $s-(|L'|-6)$ zero-sum $3$-sets we find  in $W$. Hence all the other elements in $\Gamma\setminus\{0,\iota_{\Gamma})\}$ form zero-sum $2$-subsets,  we are done.~\qed

We obtain the following:
\begin{cor} Any digraph\/ $\overrightarrow{G}$ of order\/ $n$ with no weakly connected components
 of cardinality less than\/ $4$ has a $\Gamma$-irregular labeling for every\/ $\Gamma$ such that\/
   $|\Gamma|\geq n+6$.\end{cor}
\textit{Proof.} By Corollary~\ref{bez} we can assume that $|I(\Gamma)|=1$ and $\iota_{\Gamma}\in\Gamma$ is the involution. Let $\{\overrightarrow{C}_i\}_{i=1}^t$ be the weakly connected components of $\overrightarrow{G}$. By Lemma~\ref{zerosum}, there exists a $\Gamma$-irregular labeling of a digraph $\overrightarrow{G}$ with weakly connected components $\{\overrightarrow{C}_i\}_{i=1}^t$ if and only if there exist in $\Gamma$ pairwise disjoint subsets $\{S_i\}_{i=1}^t$ such that $|S_i|=|V(\overrightarrow{C}_i)|$ and $\sum_{s\in S_i}s=0$ for every $i$, $1 \leq $i$ \leq t$. Let $m_i=|V(\overrightarrow{C}_i)|$ for every $i$, $1 \leq i \leq t$, and let $m_{t+1}=|\Gamma\setminus \{0,\iota_{\Gamma}\}|-n$. Since $m_{t+1}\geq 4$, using the sequence $\{m_i\}_{i=1}^{t+1}$, by Theorem~\ref{zero-sum}, we get the result.~\qed

\bibliographystyle{plain}

\end{document}